\documentclass[journal]{IEEEtran}
\usepackage{diagbox}
\usepackage{bbm}
\usepackage{amsfonts}
\usepackage[numbers,sort&compress]{natbib}
\usepackage{amsmath}
\usepackage[justification=centering]{caption}
\usepackage{mathrsfs}
\usepackage{algorithm}
\usepackage{algorithmic}
\usepackage{graphics}
\usepackage{multicol,multicap}
\usepackage{mathptmx}
\usepackage{times}
\usepackage{amssymb}
\usepackage{colortbl}
\usepackage{stfloats}
\usepackage{enumerate}
\usepackage{latexsym}
\usepackage{psfrag}
\usepackage{graphicx,subfigure}
\usepackage{leqno}
\usepackage{bm}
\usepackage{leftidx}
\usepackage{float}
\usepackage{blkarray}
\usepackage{arydshln}
\usepackage{setspace}
\usepackage{extarrows}
\usepackage{texsort}
\newtheorem{theorem}{\noindent{Theorem}}
\newtheorem{definition}{Definition}

\newtheorem{example}{Example}
\newtheorem{assumption}{Assumption}

\newtheorem{lemma}{Lemma}
\newtheorem{remark}{Remark}

\begin{document}

\title{A RNNs-based Algorithm for Decentralized-partial-consensus Constrained Optimization}

\author{Zicong Xia, Yang Liu, Jianlong Qiu, Qihua Ruan, and Jinde Cao
	\thanks{This work was partially supported by the National Natural Science Foundation of China under grants 11671361, 61573096, 61833005 and 61877033, the Natural Science Foundation of Zhejiang Province of China under grant LR20F030001 and D19A010003, and the Key Research and Development Project of Shandong Province of China under Grant 2017GGX10143.}
	\thanks{Z. Xia and Y. Liu are with College of Mathematics and Information Science, Zhejiang Normal University, Jinhua 321004, China (e-mail: 201531700128@zjnu.edu.cn; liuyang@zjnu.edu.cn).}
	\thanks{J. Qiu is with the School of Automation and Electrical Engineering, Key Laboratory of Complex Systems and Intelligent Computing in Universities of Shandong, Linyi University, Linyi 276005, China (e-mail: qiujianlong@lyu.edu.cn).}
	\thanks{Q. Ruan is with the School of Mathematics and Finance, Putian University, Putian 351100, China (e-mail: ruanqihua@163.com).}
	
	\thanks{J. Cao is with Jiangsu Provincial Key Laboratory of Networked Collective Intelligence, and the School of Mathematics, Southeast University, Nanjing 210096, China (e-mail: jdcao@seu.edu.cn).}}

\maketitle
\IEEEpeerreviewmaketitle
\thispagestyle{empty}
\pagestyle{empty}
\begin{abstract}
 This technical note proposes the decentralized-partial-consensus optimization with inequality constraints, and a continuous-time algorithm based on multiple interconnected recurrent neural networks (RNNs) is derived to solve the obtained optimization problems. First, the partial-consensus matrix originating from Laplacian matrix is constructed to tackle the partial-consensus constraints. In addition, using the nonsmooth analysis and Lyapunov-based technique, the convergence property about the designed algorithm is further guaranteed. Finally, the effectiveness of the obtained results is shown while several examples are presented.
\end{abstract}

\begin{IEEEkeywords}
Decentralized-partial-consensus optimization, nonsmooth analysis, continuous-time algorithms, partial-consensus matrix, recurrent neural networks (RNNs).
\end{IEEEkeywords}

\section{Introduction}

In recently decades, a great deal of the existing optimization algorithms with centralized cost functions have been developed \cite{itc1323,itac5296,is39,pieeelihuan,itnn812,itnn558}. Meanwhile, distributed optimization has also captured a major number of attention  \cite{np,sysyh,1983,ACC2019,itac10,itcns74,itcyb88490,itnnls3339,itnnls2104,is438} due to its great potential applications \cite{itac48,is136,itac1291,itc3116,itc2058}. Actually, the chief objective of distributed optimization is to design an algorithm to derive the optimal points to the optimization problem, where each agent has personal information which is only known by itself. Moreover, the methods from multi-agent theory appeal tremendous interests on account of the theoretic significance related to many kinds of fields and applications: resource assignment in communication networks \cite{itac5296,itc3116,a222}, multiagent networks \cite{itcyb88490,itnnls2104,itac48}, multi-robot motion planning, and machine learning \cite{itc2058,arc2019,YangCooperative,itcmy1}. 

As mentioned in first paragraph, the methods and theories in optimization are developed rapidly and fruitfully.   Recurrent neural networks (RNNs) are powerful tools in the handling of optimization problems. Its early application is proposed by \cite{itc1323,itcs533,itnn318,itnn558,itnn812,itnn1308}.  Many decades have witnessed the great development and improvement in RNNs for solving optimization problems. For instance, \cite{itnn1747} presents a collective neurodynamic approach with multiple interconnected RNNs for distributed optimization. \cite{itnn2344} gives a collective neurodynamic optimization method to nonnegative matrix factorization. \cite{sysyh,itcmy1} propose the RNNs-based algorithms to solve optimization and distributed optimization in quaternion field. Different from the applications above, the RNNs exchange information with their neighbors to achieve partial consensus rather than global consensus in this paper.

Various kinds of constraints have been investigated for optimization problems while considering different pratical cases in real life. For example, when it comes to the resource allocation problem, the choice of each agent localizes in a certain range, while the agents may not want to share their private information with others \cite{itac5227}. Thus, the local constraints are considered. In socialty networks, the limitations of communication capacities for each agent should also cause the constraints. In addition, Some engineering task involving construction ability, technical restrictions and time limitions have more complex constraints. In last few decades, several different constraints are considered: Inequality constraints \cite{itac1753,itcns74,itac151164}, equalities constraints \cite{is39,a222,itnn812,itac151164}, bound constraints \cite{itnn812,itac10}, local constraints \cite{itac5227,is438} and approximate constraints \cite{itcsma2957156}. In this note, we study the bound constraints and inequality constraints.

As a matter of fact, in large-scale optimization problems, the computation capacity of a single agent could not be enough to handle all the constraints of the agents, because of the performance limitations of the agents in communication capacities as well as task requirements of privacy and security. Thus, distributed optimization gets compelling attention owing to its larger capacities. Distributed optimization design are often necessary to be solved by the identical optimal point with the identical dimension for each cost function. Nevertheless, in many cases, the minimizers of different cost functions may not share the same-dimension. For instance, in a simple translation problem, the configuration of different conveyances could be different and hence the optimal transport volume for each conveyance is not identical and determined by distinct functions with respect to diverse variables. However, in this case, some variable should be identical such as the distance. The case mentioned above can be described as a decentralized-partial-consensus optimization (DPCO) problem which has the same form as cost function in distributed optimization but optimal points with different dimensions and partial consensus components. In order to be intuitionistic, Fig. 1 shows the core concept of three types of mentioned optimization.  Due to the difference for the dimensions of optimal points, the DPCO problems are more flexible and suit for more practical problems.

\begin{figure}[htp]\centering
	\includegraphics[width=9cm,height=5cm]{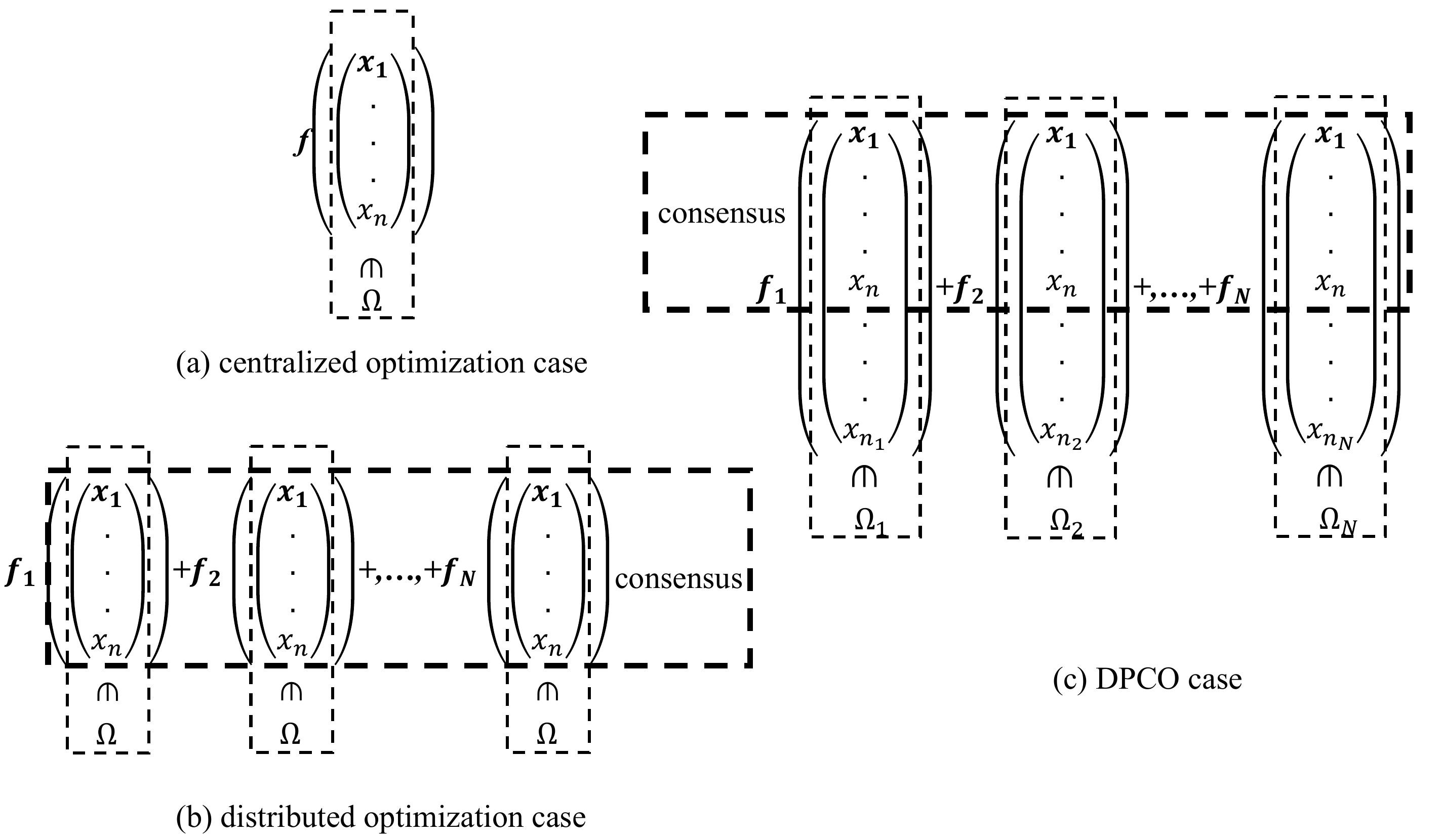}
	\caption{The distinctions among three types of optimization.}
\end{figure}
This paper aims at solving the DPCO problems by the designed continuous-time algorithm based on RNNs. Specifically, the contributions of this paper are made as summarized below.	\begin{enumerate}
	\item The model of the DPCO problems is constructed, and the partial-consensus matrix is proposed to tackle the partial-consensus constraints.
	\item In order to solve the DPCO problems, a continuous-time algorithm is designed based on multiple interconnected RNNs, and its convergence theorems are strictly proved by nonsmooth analysis and Lyapunov functions.
\end{enumerate}

The remainder of this paper is arranged as follows: Section II provides some necessary notations, definitions and lemmas including the concepts of the serial number subset and the partial-consensus matrix. In Section III, the DPCO problem is formulated, and a decentralized continuous-time algorithm is presented. Meanwhile, the description of the multiple interconnected RNNs are provided. Furthermore, the complete proof is presented to show the convergence of the designed algorithm. Section IV presents an example to demonstrate the efficiency of the results. Finally, Section V marks a brief conclusion and proposes some directions for the researches in future.
\section{Preliminaries}
\subsection{Notations}

$\mathbf{R}$ denotes the real numbers. $\mathbf{N}_{+}$ denotes the positive integers. $int(\Omega$) denotes the interior of the set $\Omega$. Define $A$ a matrix while $\delta_{min}(A) $ and $\delta_{max}(A) $ denote its smallest and largest eigenvalue respectively. $L_{N}$ denotes the $N$-dimension Laplacian matrix. $\otimes$ denotes the Kronecker product. $||\cdot||$ denotes the 2-norm. $\omega \subset \mathbf{N}_{+} $ is a finite set, $|\omega|$ denotes the cardinality of $\omega$ and $\bar{\omega}$ denotes the ordered set of $\omega$. Different from normal set, the order of ordered set can not be changed. $\bar{\omega}^{p}$ is the $p$-th entry of $\bar{\omega}$. $\bar{\omega_{i}}=\overline{\{i_{1},...,i_{n}\} }\subset \mathbf{N}_{+}$ and $\bar{\omega_{j}}=\overline{\{j_{1},...,j_{m}\} }\subset \mathbf{N}_{+}$ are two different ordered sets and all the entries are distinguished. The operator $\bar{\cup}$ leads that $\bar{\omega_{i}} \bar{\cup} \bar{\omega_{j}}=\overline{\{i_{1},...,i_{n},j_{1},...,j_{m}\}}.$ Define $col[x_{1},...,x_{n}]=[x_{1}^{T},...,x_{n}^{T}]^{T}.$

\subsection{Serial number subset and partial-consensus matrix}
$v=\{1,...,m\}$ denotes the serial number set of the components of a m-dim vector and $v_{n} \subset v$ denotes the serial number subset of the partial components of a vector. When we set $v_{n}=\{i_{1},...,i_{n}\}$ and $x=[x_{1}^{T},...,x_{m}^{T}]^{T}\in\mathbf{R}^{m}$ ($m\ge n$), then vector $x^{(v_{n})}=col[x_{i_{1}},...,x_{i_{n}}]\in \mathbf{R}^{n}$. Without loss of generality, $v_{n}$ can be supposed as $\{1,...,n\}$ and $x^{(v_{n})}=col[x_{1},...,x_{n}]\in \mathbf{R}^{n}$. Actually, the serial number subset is aimed at selecting some components to consist a new vector from a vector. First, we give a set of vectors $\{x_{i} \in \mathbf{R}^{n_{i}}|i=1,...,N\}$ and the dimensions of the vectors need not be identical. Let $\mathbf{x}=col[x_{1} ,...,x_{N} ] $ as superposition of $x_{i}(i=1,...,N)$. We set $n_{\min }=\min \left\{n_{i}\right\}$ and the corresponding serial number subsets with $n \leq n_{\min }$ are denoted as $v_{i n}$. $V=\{1,...,\sum_{i=1}^{N}n_{i}\}$ denotes the serial number set of the components of $\mathbf{x}$.
For convenience, we set $\bar{N}=\sum_{i=1}^{N}n_{i}$.

Then we propose a matrix dimension extension method by the operator
$\mathbf{E}_{\bar{\omega}}:M^{p\times p}\rightarrow M^{(p+|\omega|)\times (p+|\omega|)}$. When assume that $|\omega|=1$ and $t \in \omega$, $\mathbf{E}_{\bar{\omega}}(M)$ denotes the matrix after adding a zero row and a zero column behind the $t-1$-th row and column. When assume that $|\omega|>1$, $\mathbf{E}_{\bar{\omega}}(M)=\mathbf{E}_{\{\bar{\omega}^{|\omega|}\}}(\mathbf{E}_{\{\bar{\omega}^{|\omega|-1}\}}(....(\mathbf{E}_{\{\bar{\omega}^{1}\}}(M))))$ .

Assume that $\bar{V}_{n}=\bar{\cup}_{j=2}^{N}(\sum_{i=1}^{j-1}n_{i}+\bar{v}_{jn})\bar{\cup} \bar{v}_{1n}$ and $\overline{V\setminus V_{n}}=\bar{V}\setminus \bar{V}_{n}$.
Now, we construct the partial-consensus matrix $K_{[\Omega, n]} =\mathbf{E}_{\overline{V\setminus V_{n}}}(L_{N} \otimes I_{n})\in \mathbf{M}^{\bar{N} \times \bar{N}}(n\le n_{min})$ where $\Omega=\{x_{i}\in \mathbf{R}^{n_{i}}|i=1,...,N\}.$ For saving notations, $K_{n}$ is adopted to represent $K_{[\Omega, n]}$ and however $\Omega$ should be always considered for $K_{n}$.

\begin{lemma}\label{k1}
	Set $n_{0}=0$. When $ K_{n(i,j)}\ne 0$ if and only if $\exists N_{1},N_{2} \in \mathbf{N}_{+}$ and $N_{1}\le N,N_{2}\le N$ such that $i-j=\sum_{t=0}^{N_{1}-1}n_{t}-\sum_{t=0}^{N_{2}-1}n_{t}$.
\end{lemma}
\begin{IEEEproof}
	It is noticed that $L_{n} \otimes I_{N(i,j)}\ne0$ if and only if $\exists N_{3} \in \mathbf{N}_{+}$ and $|N_{3}| \le N-1$ such that $i-j=N_{3}n$. \\(Necessity) According to the operator $\mathbf{E}$, for $\sum_{i=t}^{N_{1}-1}n_{t}\le i\le\sum_{t=1}^{N_{1}}n_{t}\ $ and $\sum_{i=t}^{N_{2}-1}n_{t}\le j\le\sum_{t=1}^{N_{2}}n_{t}\ $, $K_{n(i,j)}=L_{N} \otimes I_{n(i_{0},j_{0})}$, where $i_{0}=i-\sum_{t=1}^{N_{1}-1}(n_{t}-n)$ and $j_{0}=j-\sum_{t=1}^{N_{2}-1}(n_{t}-n)$.
	The providing the range of $i_{0},j_{0}$: $(N_{1}-1)n\le i_{0}\le N_{1}n$ and $(N_{2}-1)n\le j_{0}\le N_{2}n$. From the condition in which $K_{n(i,j)}=L_{N} \otimes I_{n(i_{0},j_{0})}\ne0$ , there exists $N_{3}$ such that $i_{0}-j_{0}=N_{3}n$. Then considering about the range of $i_{0},j_{0}$, one can find that $N_{1}=N_{2}+N_{3}$. Next, $i_{0}-j_{0}=i-\sum_{t=1}^{N_{1}-1}(n_{t}-n)-(j-\sum_{t=1}^{N_{2}-1}(n_{t}-n))=N_{3}n$, combining with $N_{1}=N_{2}+N_{3}$, we have  $i-j=\sum_{t=1}^{N_{1}-1}n_{t}-\sum_{t=1}^{N_{2}-1}n_{t}=\sum_{t=0}^{N_{1}-1}n_{t}-\sum_{t=0}^{N_{2}-1}n_{t}$.
	\\(Sufficiency) Suppose that $\exists N_{1},N_{2} \in N_{+}$ and $N_{1}\le N,N_{2}\le N$ such that $i-j=\sum_{t=0}^{N_{1}-1}n_{t}-\sum_{t=0}^{N_{2}-1}n_{t}$. It is obvious that $\exists N_{4},N_{5} \in N_{+}$ and $N_{1}\le N,N_{2}\le N$, such that $\sum_{i=t}^{N_{4}-1}n_{t}\le i\le\sum_{t=1}^{N_{4}}n_{t}\ $, $\sum_{i=t}^{N_{5}-1}n_{t}\le j\le\sum_{t=1}^{N_{5}}n_{t}\ $, and $i-j=\sum_{t=0}^{N_{1}-1}n_{t}-\sum_{t=0}^{N_{2}-1}n_{t}=\sum_{t=0}^{N_{4}-1}n_{t}-\sum_{t=0}^{N_{5}-1}n_{t}$. In fact, when selecting $i_{0}=i-\sum_{t=1}^{N_{4}-1}(n_{t}-n)$ and $j_{0}=j-\sum_{t=1}^{N_{5}-1}(n_{t}-n)$, they lead to that $K_{n(i,j)}=L_{N} \otimes I_{n(i_{0},j_{0})}$. Next, we calculate $i_{0}-j_{0}=i-\sum_{t=1}^{N_{4}-1}(n_{t}-n)-(j-\sum_{t=1}^{N_{5}-1}(n_{t}-n))=i-j-(\sum_{t=0}^{N_{1}-1}n_{t}-\sum_{t=0}^{N_{2}-1}n_{t})+(N_{4}-N_{5})n=(N_{4}-N_{5})n$. Finally, set $N_{3}=N_{4}-N_{5}$, then $K_{n(i,j)}=L_{N} \otimes I_{n(i_{0},j_{0})}\ne0$ holds.
\end{IEEEproof}
\begin{lemma}\label{k2}
$K_{n}\mathbf{x}=0$ if and only if $x_{i}^{(v_{n})}=x_{j}^{(v_{n})}$, $i,j=1,...,N.$
\end{lemma}
\begin{IEEEproof}
	To commence with the equation $K_{n}\mathbf{x}=0$, it can be easily derived that $\sum_{j=1} ^{\bar{N}}K_{n(i,j)}\mathbf{x}_{j}=0$ for $i=1,...,N$. According to the properties of $L_{N} \otimes I_{n}$, we find that $K_{n(i,i)}=-\sum_{i=1}^{\bar{N}}K_{n(i,j)}$. On account of Lemma \ref{k1}, considering the case in which $i= \bar{V}_{n}^{1} $, we have
	\begin{equation}\label{kk}
	\begin{array}{l}K_{n(\bar{V}_{n}^{1},\bar{V}_{n}^{1})}\mathbf{x}_{\bar{V}_{n}^{1}}+K_{n(\bar{V}_{n}^{1},n_{1}+\bar{V}_{n}^{1})}\mathbf{x}_{\bar{V}_{n}^{1+{n}_{1}}}+...\\+K_{n(\bar{V}_{n}^{1},\sum_{t=1}^{N-1}n_{t}+\bar{V}_{n}^{1})}\mathbf{x}_{\bar{V}_{n}^{1+\sum_{t=1}^{N-1}n_{t}}}\\=\sum_{j=2}^{N}K_{n(\bar{V}_{n}^{1},\bar{V}_{n}^{1+\sum_{t=1}^{j-1}n_{t}})}(\mathbf{x}_{\bar{V}_{n}^{1+\sum_{t=1}^{j-1}n_{t}}}-\mathbf{x}_{\bar{V}_{n}^{1}})=0.
	\end{array}
	\end{equation}
	Then choose $i=\bar{V}_{n}^{1+n_{1}},...,\bar{V}_{n}^{1+\sum_{t=1}^{j-1}n_{t}}$, respectively, it leads the set of a equations from equation (\ref{kk}):
	
	\begin{equation}\label{kkk}
	\left\{\begin{array}{l}
	\sum_{j\ne 1} K_{n(\bar{V}_{n}^{1},\bar{V}_{n}^{1+\sum_{t=1}^{j-1}n_{t}})}(\mathbf{x}_{\bar{V}_{n}^{1+\sum_{t=1}^{j-1}n_{t}}}-\mathbf{x}_{\bar{V}_{n}^{1}})=0,\\\sum_{j\ne 2} K_{n(\bar{V}_{n}^{1+n_{1}},\bar{V}_{n}^{1+\sum_{t=1}^{j-1}n_{t}})}(\mathbf{x}_{\bar{V}_{n}^{1+\sum_{t=1}^{j-1}n_{t}}}-\mathbf{x}_{\bar{V}_{n}^{1+n_{1}}})=0,\\~~~~~~~~~~~~~~~~~~.\\~~~~~~~~~~~~~~~~~~.\\~~~~~~~~~~~~~~~~~~.\\\sum_{j\ne N} K_{n(\bar{V}_{n}^{1+\sum_{t=1}^{N-1}n_{t}},\bar{V}_{n}^{1+\sum_{t=1}^{j-1}n_{t}})}(\mathbf{x}_{\bar{V}_{n}^{1+\sum_{t=1}^{j-1}n_{t}}}-\mathbf{x}_{\bar{V}_{n}^{1+\sum_{t=1}^{N-1}n_{t}}})=0.\end{array}\right.
	\end{equation}
	
	It is easy to observe that (\ref{kkk}) guarantees $\mathbf{x}_{\bar{V}_{n}^{1}}=...=\mathbf{x}_{\bar{V}_{n}^{1+\sum_{t=1}^{j-1}n_{t}}}$. Due to the construction of $\bar{V}_{n}$ and with the further derivation, one can find $\bar{V}_{n}^{1+\sum_{t=1}^{j-1}n_{t}}=\sum_{t=1}^{j-1}n_{t}+\bar{V}_{n}^{1}=\sum_{t=1}^{j-1}n_{t}+\bar{v}_{1n}^{1}$, then $\mathbf{x}_{\bar{V}_{n}^{1+\sum_{t=1}^{j-1}n_{t}}}=x_{j}^{1}=x_{j}^{(v_{n})^{1}}$ holds. Finally, it leads to the result:	${x_{1}^{(v_{n})}}^{1}=...={x_{N}^{(v_{n})}}^{1}$. Similarly, ${x_{1}^{(v_{n})}}^{p}=...={x_{N}^{(v_{n})}}^{p}$ for $p=1,...,n$. This completes the conclusion: $x_{i}^{(v_{n})}=x_{j}^{(v_{n})}$.

\end{IEEEproof}

From Lemma \ref{k2}, we obtain that the main function of $K_{n}$ is to make some parts of several vectors identical and then its properties originating from the laplacian matrix are proved as follows:
\begin{lemma}
	Some properties of $K_{n}$ are mentioned:
	\begin{enumerate} [1)]
		
		\item $K_{n}$ is positive semi-definite.
		
		\item $K_{n}$ has $nN$ nonnegative eigenvalues. $\delta_{min}(K_{n})=\delta_{min}(L_{N} \otimes I_{n})=0$ and its multiplicity is $\bar{N}-nN+1$. $\delta_{max}(K_{n})=\delta_{max}(L_{N} \otimes I_{n})$.

	\end{enumerate}\end{lemma}	\begin{IEEEproof}
	The conclusion can be easily proved by the characteristic polynomial of $K_{n}$ which is similar to the one of $L_{N} \otimes I_{n}$. Thus, the proof is omitted.
\end{IEEEproof}

Now, we give a simple example to explain the construction of the partial-consensus matrix from a group of vectors.
\begin{example}
	We give a set of vectors
	$$\begin{array}{l}
	x_{1}=\left\{x_{11}, x_{12}, x_{13}\right\}, x_{2}=\left\{x_{21}, x_{22}, x_{23}, x_{24}\right\},\\
	x_{3}=\left\{x_{31}, x_{32}, x_{33}, x_{34}, x_{35}\right\}.
	\end{array}$$
	
	Define $\Omega=\{x_{1},x_{2},x_{3}\}$. It is easy to get $n_{min}=3$ and set $n=3$ which leads to $v_{3}=\{1,2,3\}.$ $\mathbf{x}$ can be also obtained by superposing $x_{1},x_{2},x_{3}$ and $V=\{1,...,12\}.$ Then $V_{3}=\{1,2,3,4,5,6,8,9,10\}$ and 
	$$\mathbf{x}^{(V_{3})}=\{x_{11},x_{12},x_{13},x_{21},x_{22},x_{23},x_{31},x_{32},x_{33}\}, $$
	$$\mathbf{x}^{(V\setminus V_{3})}=\{x_{24},x_{34},x_{35}\}. $$ We select the Laplacian matrix $L_{3}=\left[
	\begin{array}{ccc}
	2& -1 & -1\\
	-1& 2 & -1\\
	-1 & -1 & 2\\
	\end{array}
	\right]  $. Then $K_{3}$ can be selected as follows:
	\footnotesize{	$$
		K_{3}=\left[
		\begin{array}{cccccccccccc}
		2 & 0 & 0& -1& 0& 0& 0& -1& 0& 0& 0& 0\\
		0& 2 & 0& 0& -1& 0& 0& 0& -1& 0& 0& 0\\
		0 & 0 & 2 & 0& 0& -1& 0& 0& 0& -1& 0& 0\\
		-1& 0 & 0& 2& 0& 0& 0& -1& 0& 0& 0& 0\\
		0& -1 & 0& 0& 2& 0& 0& 0& -1& 0& 0& 0\\
		0 & 0 & -1& 0& 0& 2& 0& 0& 0& -1& 0& 0\\
		0& 0 & 0& 0& 0& 0& 0& 0& 0& 0& 0& 0\\
		-1& 0 & 0& -1& 0& 0& 0& 2& 0& 0& 0& 0\\
		0& -1& 0& 0& -1& 0& 0& 0& 2& 0& 0& 0\\
		0& 0 & -1& 0& 0& -1& 0& 0& 0& 2& 0& 0\\
		0 & 0 & 0& 0& 0& 0& 0& 0& 0& 0& 0& 0\\
		0 & 0 & 0& 0& 0& 0& 0& 0& 0& 0& 0& 0\\
		\end{array}
		\right] .
		$$}
	\normalsize{The graph of $L_{3}$ and $K_{3}$'s corresponding adjacency matrix is shown in Fig. \ref{to1}.}

	 \begin{figure}[!hbt]\centering
	 	\includegraphics[width=6cm,height=4.2cm]{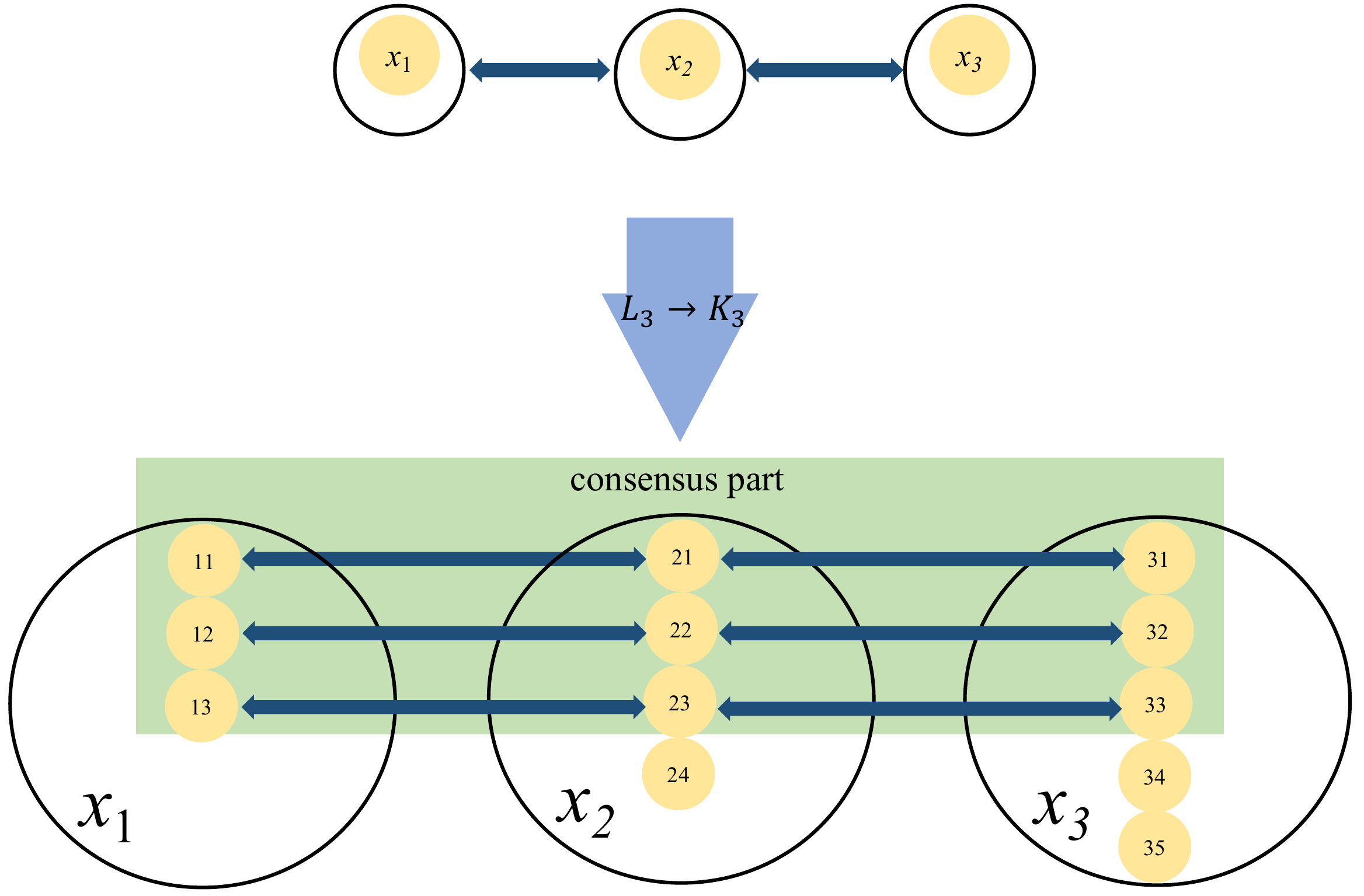}
	 	\caption{The graph of corresponding adjacency matrix of $L_{3}$ and $K_{3}$ in Example 1.}
	 	\label{to1}
	 \end{figure}

\end{example}

In Section II, let $v_{n} \subset v$ denote the serial number subset of the partial components of a vector. Correspondingly, for a set of vectors $\Omega=\{x_{1},x_{2},...,x_{N}\}$, the consensus parts are $x_{1j}=x_{2j}=,...,=x_{Nj}, j=1,...,n.$ In essence, we can also choose any n components of each vector $\in \Omega$ and just move the selected components to the top n positions of each vector $\in \Omega$. To achieve the mentioned movement, we adopt the more general form of $v_{n}=\{i_{1},...,i_{n}\}$ and propose the definition on the permutation matrix $A[v\setminus v_{n}]$ which can change $x \in \mathbf{R}^{m}$ into $col\{x^{(v_{n})},x^{(v\setminus v_{n})}\}$. Moreover, in the remainder of this note, we just the case in $v=\{1,...,m\}$. Thus, the permutation matrix $A[v\setminus v_{n}]$ is not used, but it really plays the role of random selection of consensus parts.
\begin{lemma}\label{l}
	$col\{x^{(v_{n})},x^{(v\setminus v_{n})}\}=A[v\setminus v_{n}]x$,  where	
	$$A[v\setminus v_{n}]_{(ij)}=\left\{\begin{array}{ll}
	1 & (i,j)=(\bar{v}^{p},\bar{\omega}^{p}),p=1,...,m, \\0 &  \text {others,}
	\end{array}\right.$$
	with $\bar{\omega}=\bar{v}_{n}\bar{\cup}\overline{v\setminus v_{n}}.$
\end{lemma}
\begin{IEEEproof}
	In fact, the construction of the permutation matrix is similarly to the normal permutation matrix, the proof is easy and hence is omitted.
\end{IEEEproof}
\subsection{Projection and nonsmooth analysis in convex optimization}
Suppose that $\Omega$ is a convex set, and $P_{\Omega}(x)$ denotes the projection of x onto $\Omega$, i.e., $P_{\Omega}(x)=argmin_{x^{\prime}\in \Omega}||x-x^{\prime}||$. Then we have the property \cite{np}:
\begin{equation}\label{pi}
\left(P_{\Omega}(x)-x\right)^{T}\left(P_{\Omega}(x)-x^{\prime}\right) \leq 0, \forall x \in \mathbb{R}^{n}, \forall x^{\prime} \in \Omega.
\end{equation}
From \cite{np}, $\operatorname{cone}(\Omega)=\{\gamma x: x \in \Omega, \gamma \geq 0\}$ denotes the convex cone of $\Omega$. $N_{\Omega}(x)=[\operatorname{cone}(\Omega-x)]^{\circ}$ is the normal cone of $\Omega$ at $x \in \Omega$.
\begin{lemma} Let $X$ be a closed convex set and let $x \in \Omega .$ Then
	\begin{equation}\label{cN} N_{\Omega}(x)=\left\{x^{\prime} \in \mathbb{R}^{n}: P_{\Omega}(x+x^{\prime} )=x\right\}.\end{equation}
\end{lemma}
\begin{definition} (generalized gradient) A function $f: \mathbb{R}^{n} \rightarrow \mathbb{R}$ is Lipschitz continuous in $X . v$ is a vector in $X, \partial f(x)=$ $\left\{\xi \in \mathbb{R}^{n} \mid \xi \in X: f^{0}(x ; v) \geq\langle\xi, v\rangle\right\}$ denotes the generalized
	gradient of $f$ at $x,$ where $f^{0}(x ; v)$ is finite and well defined, and $f^{0}(x ; v)=\lim _{x \rightarrow 0} \sup ((f(x+t v)-f(x)) / t)$.
\end{definition}

$\partial_{\mu} f(x)$($\mu$ is partial of $x$) denotes the subgradient of the partial derivatives, specificly defined as follows,
$$
\partial_{\mu} f(x)=\left\{\xi_{1} \in \mathbb{R}^{n}: f(\mu, y)-f(v, y) \leq\left\langle\xi_{1}, \mu-v\right\rangle\right\}.
$$
\begin{lemma}
	If $F: \mathbb{R}^{n} \rightarrow \mathbb{R}$ is $L$-Lipschitz continuous, then $\left\{x_{i}\right\}$ be a sequence in $\mathbb{R}^{n}$ and $\xi_{i} \in \partial F\left(x_{i}\right).$ Suppose $x_{i}$ converges to $x$ and $\xi$ is a cluster point of $\xi_{i},$ hence $\xi \in \partial F(x)$.
\end{lemma}
\begin{remark}
	The concepts related to nonsmooth analysis based on generalized gradient and differential inclusion are analyzed in \cite{np,itac5227}. The generalized gradient with respect to $x$ is often a set of vectors whose dimension is the same as $x$. When it comes to the inequality scaling, it is usual to choose a feasible element from the generalized gradient like formula (3a) in \cite{itac5227}. In this paper, for convenience and saving notation, when the inequalities involve the generalized gradient, the process of inequality scaling in the proof just continues without selecting some feasible elements from the generalized gradient.
\end{remark}
\section{Main Results}
In this Section, we propose the concept of DPCO and give a RNNs-based algorithm to solve DPCO problems. Compared with the distributed optimization, the cost functions of DPCO has the same form as ones in distributed optimization but its optimal points have different dimensions and partial consensus components. Specifically, compared with the distributed optimization and decentralized optimization, there are two main differences in the DPCO:
\begin{enumerate}
	\item The dimensions of the optimal point for each objective function need not be identical.
	\item Only some components of the optimal points for each objective function are identical.
\end{enumerate}

To this end, this note mainly considers the decentralized-partial-consensus optimization problem modeled as follows:
\begin{equation} \label{jbop}
\begin{aligned}\min F(\mathbf{x}) =&\sum_{i=1}^{N}f_{i}(x_{i}),\\
~s.t.~&x_{i}^{(v_{n})}=x_{j}^{(v_{n})},\\&g_{i}(x_{i})\le 0, \\& x_{i} \in \Omega_{i}, i,j=1,...,N.
\end{aligned}
\end{equation}
Where $x_{i} \in \mathbf{R}^{n_{i}}$, $\mathbf{x}=col[x_{1} ,...,x_{N} ] $, $ n \le min\{n_{i}\}$, $f_{i}: \mathbf{R}^{n_{i}}\rightarrow \mathbf{R}$ and $g_{i}: \mathbf{R}^{n_{i}}\rightarrow \mathbf{R}^{m_{i}}$, $i=1,...,N.$

In this note, we call formula (\ref{jbop}) as Problem (\ref{jbop}).
Actually, it can be obtained that Problem (\ref{jbop}) can be changed into distributed optimization problem, when $n_{i}=n_{j}, i,j=1,...,N$, $v_{n}=v$ and $m_{i}=m_{j}, i,j=1,...,N$. Similarly, Problem (\ref{jbop}) can be also changed into decentralized optimization. Thus, the considered problems are more general than distributed optimization and decentralized optimization problems.

For brevity and $\forall i=1,...,N$, denote $ \lambda=col[\lambda_{1}, \ldots,\lambda_{N}]$ with $\lambda_{i} \in \mathbf{R}^{n_{i}}$, $ \mu=col[\mu_{1}, \ldots, \mu_{N}]$ with $\mu_{i} \in \mathbf{R}^{m_{i}}$, $P_{+}^{\mu}=col[P_{+}^{\mu_{1}}, \ldots, P_{+}^{\mu_{N}}]$ with $P_{+}^{\mu_{i}}=P_{\mathbb{R}_{+}^{m}}(\mu_{i}+g_{i}(x_{i}))
$, $G(\mathbf{x})=col[g_{1}\left(x_{1}\right), \ldots, g_{N}\left(x_{N}\right)] $, and $\Omega=\prod_{i=1}^{N}\Omega_{i}=\Omega_{1}\times \Omega_{2}\times ,...,\times  \Omega_{N}$. $\times$ denotes the Cartesian product.

With the compact notations, a decentralized continuous-time algorithm based on RNNs, in order to tackle the Problem (\ref{jbop}), is proposed as follows:

\begin{equation} \label{maina1}
\left\{\begin{aligned}
\dot{x}\in &2 \delta  P_{\Omega}(\mathbf{x}- \partial F(\mathbf{x})-\partial G(\mathbf{x}) P_{+}^{\mu}-K_{n}\lambda \\&-K_{n}\mathbf{x} )-2\delta \mathbf{x},
\\\dot{\lambda}=&K_{n}\mathbf{x} ,\\
\dot{\mu}=&\delta \left(P_{+}^{\mu}-\mu\right) ,
\end{aligned}\right.
\end{equation}
where $\delta=\delta_{max}(K_{n})+1$.

Then turn the Algorithm (\ref{maina1}) in compact form into a coupling form to describe the RNNs clearly:
\begin{equation}
\left\{\begin{aligned}
\frac{d x_{i}}{d t} \in & 2\delta[P_{\Omega_{i}}(x_{i}-\partial f_{i}(x_{i})-\partial g_{i}(x_{i})P_{+}^{\mu_{i}}\\&-\sum_{j=1, j \neq i}^{\bar{N}} K_{n(i,j)}(x_{i}-x_{j}+\lambda_{i}-\lambda_{j})-x_{i}],\\
\frac{d z_{i}}{d t}=&\sum_{j=1, j \neq i}^{\bar{N}} K_{n(i,j)}(x_{i}-x_{j}),\\
\frac{d \mu_{i}}{d t}=&\delta(P_{+}^{\mu_{i}}-\mu_{i}). &
\end{aligned}\right.
\end{equation}

\begin{figure}[htp]\centering
	\includegraphics[width=10.4cm,height=8cm]{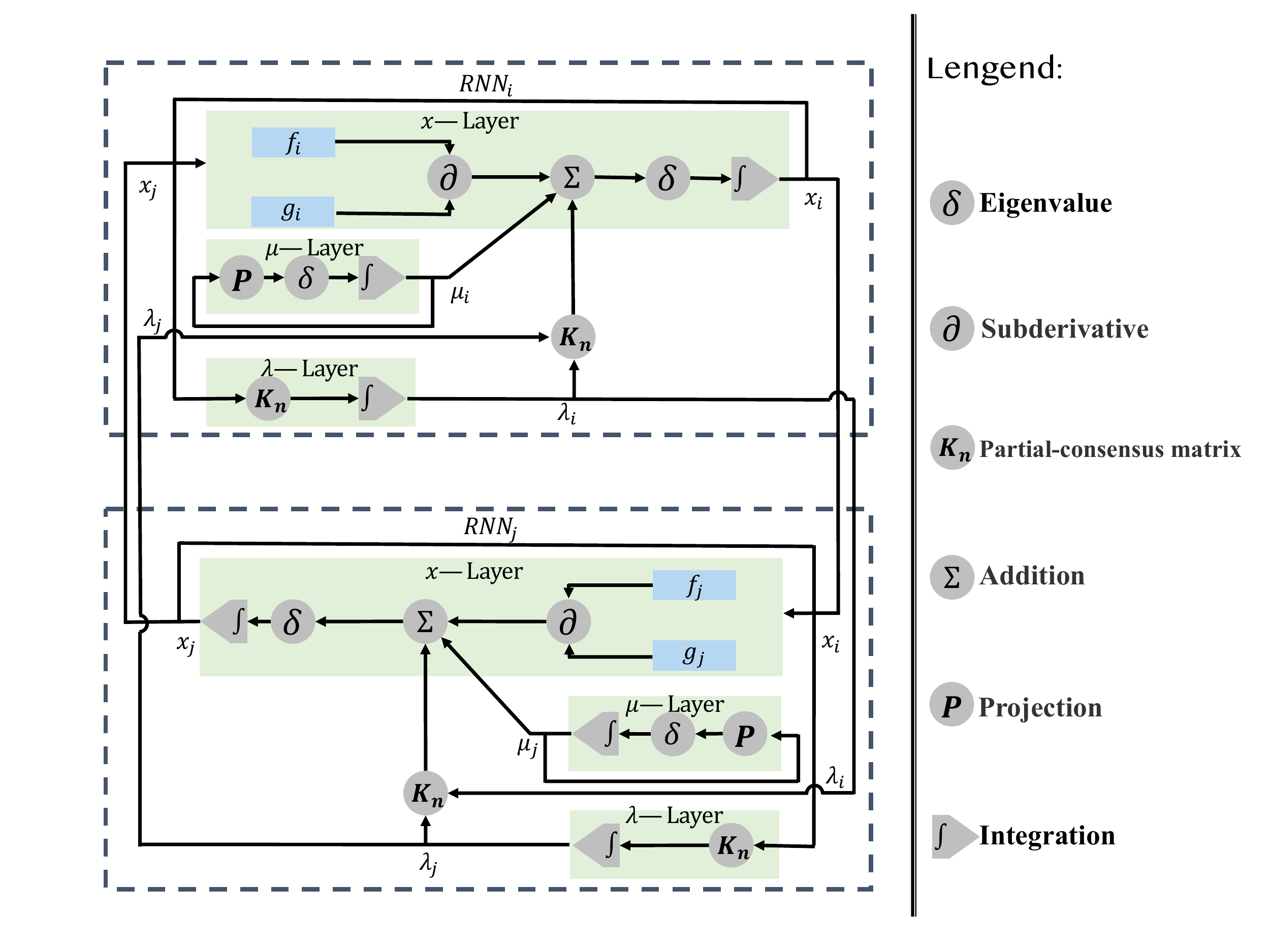}
	\caption{Link diagram of RNNs in Algorithm (\ref{maina1}).}
	\label{BP}
\end{figure}

Meanwhlie, the link diagram of the obtained multiple interconnected RNNs in Algorithm (\ref{maina1}) is shown in Fig. \ref{BP}, in which one can find that each RNN has a two-layer structure and the RNNs, for example, $i$ and $j$ are linked by the constructed partial-consensus matrix $K_{n}$. Each RNN$i$ exchanges the information on $x_{i}$ and $\lambda_{i}$ with its neighbors determined by $K_{n}$. Specifically, in the collective RNNs, $f_{i}$, $g_{i}$ and $\Omega_{i}$ in (\ref{maina1}) are only known by RNN$_{i}$, i.e., each RNN$_{i}$ is arranged to minimize a local objective function subject to its local constraints rather than global information. Different from the applications in RNNs such as \cite{itnn1747,itnn2344,itcmy1}, the RNN$_{i}$ exchanges information with its neighbors to achieve partial consensus rather than global consensus.

In the remainder of Section III, the functions for Problem (\ref{jbop}) and convergence of Algorithm (\ref{maina1}) will be illustrated. Before that, the necessary assumptions are given, which is the well-known convex optimization items to the convex inequality constraints.


\begin{assumption}\label{a1} \quad
	\begin{enumerate}[1)]	\item (Convexity and continuity) $\Omega_{i}$ is convex set, $f_{i}, g_{i}$ are convex and continunous, $i=1,...,N$.
		\item (Slater's condition) There exists $\bar{\mathbf{x}} \in int(\Omega)$ such that $g(\bar{\mathbf{x}})<0$.
	\end{enumerate}
\end{assumption}

Assumption \ref{a1} is a normal term in optimization problem. Convexity and continuity keep that the optimal points are the global solutions. Slater's condition guarantees the constraint qualification for optimization problems. In addition, the Problem (\ref{jbop}) should have at least one solution. Now, it is necessary to propose some preliminary lemmas under Assumption \ref{a1}:
\vspace{5mm}
\begin{lemma} \label{hhhha}
	Under Assumption \ref{a1}, a point $\mathbf{x}^{*}$ is the solution to Problem (\ref{jbop}) if and only if there exist $\lambda^{*}$ and $\mu^{*}$ such that
	
	\begin{equation}\label{keys}
	\begin{aligned}\mathbf{0} \in &\partial F(\mathbf{x}^{*})+\partial G(\mathbf{x}^{*}) \mu^{*}+K_{n}\lambda^{*}+N_{\Omega}(\mathbf{x}^{*}),\\&g_{i}\left(x_{i}^{*}\right) \in N_{\mathbb{R}_{+}^{m_{i}}}(\mu_{i}^{*}), \forall i=1,...,N.
	\end{aligned}
	\end{equation}
\end{lemma}
\begin{IEEEproof}
	According to Lemma \ref{k2}, it is easy to verify that Problem (\ref{jbop}) is equivalent to
	\begin{equation} \label{jbop2}
	\begin{aligned}\min F(\mathbf{x}) =&\sum_{i=1}^{N}f_{i}(x_{i}),\\
	~s.t.~K_{n}\mathbf{x}=0, g_{i}(x_{i})\le& 0, x_{i} \in \Omega_{i},i=1,...,N.
	\end{aligned}
	\end{equation}
	Then inspired by \cite{np} under Assumption \ref{a1}, there exists an optimal point $\mathbf{x}^{*}$, $ \lambda^{*}=col[\lambda_{1}^{*}, \ldots,\lambda_{N}^{*}]$ with $\lambda_{i} ^{*}\in \mathbf{R}^{n_{i}}$, and $ \mu^{*}=col[\mu_{1}^{*}, \ldots, \mu_{N}^{*}]$ with $\mu_{i}^{*} \in \mathbf{R}^{m_{i}}$ such that
	\begin{equation}\label{sla}
	\mathbf{0} \in \partial F(\mathbf{x}^{*})+\partial G(\mathbf{x}^{*}) \mu^{*}+K_{n}\lambda^{*}+N_{\Omega}(\mathbf{x}^{*}).
	\end{equation}
	\begin{equation}\label{sla1}
	K_{n}\mathbf{x}^{*}=0.
	\end{equation}
	\begin{equation}\label{sla2}
	\mu_{i}^{*} \geq 0, g_{i}\left(x_{i}\right) \leq 0,\left\langle\mu_{i}^{*},g_{i}\left(x_{i}\right)\right\rangle=0, \forall i=1,...,N.
	\end{equation}
	It is no doubt that (\ref{sla}) amounts to the first formula of (\ref{keys}), and while (\ref{sla2}) amounts to the second formula of (\ref{keys}), which keeps this proof.
\end{IEEEproof}

\vspace{5mm}

\begin{theorem}\label{e}
	Assume that Assumption \ref{a1} holds. the equilibrium point of Algorithm (\ref{maina1}) is the solution to Problem (\ref{jbop}).
\end{theorem}

\vspace{5mm}
\begin{IEEEproof}
	Suppose that $(\hat{\mathbf{x}},\hat{\lambda},\hat{\mu})$ is the equilibrium point of Algorithm (\ref{maina1}), which indicates
	\begin{subequations} \label{zong}
		\begin{equation} \label{zong1}
		\begin{aligned}
		\mathbf{0} &\in 2 \delta  P_{\Omega}(\hat{\mathbf{x}}- \partial F(\hat{\mathbf{x}})-\partial G(\hat{\mathbf{x}}) P_{+}^{\hat{\mu}}-K_{n}\hat{\lambda} -K_{n}\hat{\mathbf{x}} )-2\delta \hat{\mathbf{x}}, \end{aligned}
		\end{equation}
		\begin{equation}  \label{zong2}
		\mathbf{0}=K_{n}\hat{\mathbf{x}},
		\end{equation}
		\begin{equation}  \label{zong3}
		\mathbf{0}=\delta \left(P_{+}^{\hat{\mu}}-\hat{\mu}\right).
		\end{equation}
	\end{subequations}
	
	Thus, (\ref{zong3}) implies that $P_{+}^{\hat{\mu}}=\hat{\mu}$ which guarantees $P_{\mathbb{R}_{+}^{m}}(\hat{\mu}_{i}+g_{i}(x_{i}))=\hat{\mu}_{i}$. Then, the second formula of (\ref{keys}) holds. To sum up with $P_{+}^{\hat{\mu}}=\hat{\mu}$, (\ref{zong1}) and (\ref{zong2}), the first formula of (\ref{keys}) holds. On account of Lemma \ref{keys}, $\hat{\mathbf{x}}$ is a solution to Problem (\ref{jbop}).
\end{IEEEproof}

\vspace{3mm}
Lemma \ref{hhhha} and Theorem \ref{e} connect the equilibrium point of Algorithm (\ref{maina1}) with the solution to Problem (\ref{jbop}) by the tool in KKT (Karush-Kuhn-Tucker) condition and equilibrium point anlysis. Now, the converge of the designed algorithm can be proved by the approaches of nonsmooth analysis and Lyapunov-based technique. The proof consists of two parts: 1) The boundness of $\mathbf{x}, \lambda, \mu$ from Algorithm (\ref{maina1}); 2) Algorithm (\ref{maina1}) converge to the equilibrium point of Algorithm (\ref{maina1}). The proving method is to construct Lyapunov functions and use the properties of coercive functions.
\vspace{3mm}

\begin{theorem}\label{TT2}
	Assume that Assumption \ref{a1} holds, then the following conclusions are right:
	\begin{enumerate}
		\item $\mathbf{x}, \lambda, \mu$ from Algorithm (\ref{maina1}) are bounded;
		\item The states $\mathbf{x}$ of Algorithm (\ref{maina1}) will converge to a solution for Problem (\ref{jbop}) from any initial point $\left(\mathbf{x}(0), \lambda(0), \mu(0)\right)$.
	\end{enumerate}
\end{theorem}
\begin{IEEEproof}
	Similarly to the proof of Theorem \ref{e}, $(\hat{\mathbf{x}},\hat{\lambda},\hat{\mu})$ denotes the equilibrium of of Algorithm (\ref{maina1}). The proof begins with construsting the Lyapunov functions:
	\begin{equation} \label{11}\begin{aligned}
	&V=V_{1}+V_{2}+V_{3}+V_{4}, \\[2mm]&V_{1} =F(\mathbf{x})+\frac{1}{2}\left\|P_{+}^{\mu}\right\|^{2}+\frac{1}{2}(\mathbf{x}+\lambda)^{T}K_{n}(\mathbf{x}+\lambda),
	\\[2mm]&V_{2} =V_{1}(\mathbf{x},\lambda,\mu)-V_{1}(\hat{\mathbf{x}},\hat{\lambda},\hat{\mu})-\partial_{\mathbf{x}} V_{1}(\hat{\mathbf{x}},\hat{\lambda},\hat{\mu})^{T}(\mathbf{x}-\hat{\mathbf{x}})  \\[2mm]&~~~~~~-(\hat{\mathbf{x}}+\hat{\lambda})^{T}K_{n}(\hat{\mathbf{x}}+\hat{\lambda})-\hat{\mu}^{T}(\mu-\hat{\mu}),
	\\[2mm]&V_{3} =\frac{1}{2}\|\mathbf{x}-\hat{\mathbf{x}}\|^{2}+\frac{1}{2}(\lambda-\hat{\lambda})^{T} (2 \delta I _{\bar{N}}-K_{n}) (\lambda-\hat{\lambda} ), \\[2mm]&V_{4} =\frac{1}{2}\|\mu-\hat{\mu}\|^{2}.\end{aligned}\end{equation}
	
	Then derivative of $V_{2}$ is calculated:
	
	\begin{equation}  \begin{aligned}
	\dot{V}_{2}=&(\partial_{\mathbf{x}}V_{1}(\mathbf{x},\lambda,\mu)-\partial_{\mathbf{x}}V_{1}(\hat{\mathbf{x}},\hat{\lambda},\hat{\mu}))\dot{\mathbf{x}}\\&+(\mathbf{x}-\hat{\mathbf{x}})^{T} K_{n}\dot{\lambda}+(\lambda-\hat{\lambda})^{T} K_{n}\dot{\lambda}\\&+\left(P_{+}^{\mu}-\hat{\mu}\right)^{T} \dot{\mu}.	 \end{aligned}\end{equation}
	
	Divide $\dot{V}_{2}$ into three parts $\dot{V}_{21}$, $\dot{V}_{22}$ and $\dot{V}_{23}$:
	\begin{equation}  \dot{V}_{21}=(\partial_{\mathbf{x}}V_{1}(\mathbf{x},\lambda,\mu)-\partial_{\mathbf{x}}V_{1}(\hat{\mathbf{x}},\hat{\lambda},\hat{\mu}))\dot{\mathbf{x}},
	\end{equation}
	\begin{equation} \dot{V}_{22}=(\mathbf{x}-\hat{\mathbf{x}})^{T} K_{n}\dot{\lambda}+(\lambda-\hat{\lambda})^{T} K_{n}\dot{\lambda},	\end{equation}  \begin{equation}
	\dot{V}_{23}=\left(P_{+}^{\mu}-\hat{\mu}\right)^{T} \dot{\mu}. \end{equation}
	On the purpose of convenience, set $$\partial_{\mathbf{x}}V_{1}(\mathbf{x},\lambda,\mu)=\partial_{\mathbf{x}}V_{1},$$ $$\partial_{\mathbf{x}}V_{1}(\hat{\mathbf{x}},\hat{\lambda},\hat{\mu})=\partial_{\mathbf{x}}\hat{V}_{1},$$ $$P_{\Omega}=P_{\Omega}(\mathbf{x}- \partial F(\mathbf{x})-\partial G(\mathbf{x}) P_{+}^{\mu}-K_{n}\lambda-K_{n}\mathbf{x} ).$$
	At this moment, one can calculate that
	\begin{equation}  \begin{aligned}
	\dot{V}_{3}+\dot{V}_{21}& =2 \delta(\mathbf{x}-\hat{\mathbf{x}})^{T}(P_{\Omega}-\mathbf{x})+2 \delta(\partial_{\mathbf{x}}V_{1}-\partial_{\mathbf{x}}\hat{V}_{1})^{T}(P_{\Omega}-\mathbf{x}) \\&~~+(\lambda-\hat{\lambda})^{T}(2 \delta I _{\bar{N}}-K_{n}) \dot{\lambda}\\&=2 \delta(\partial_{\mathbf{x}}V_{1}-\partial_{\mathbf{x}}\hat{V}_{1})^{T}(P_{\Omega}-\hat{\mathbf{x}})\\&~~+2\delta(\partial_{\mathbf{x}}V_{1}-\partial_{\mathbf{x}}\hat{V}_{1})^{T}(\hat{\mathbf{x}}-\mathbf{x})
	\\&~~+2 \delta(\mathbf{x}-\hat{\mathbf{x}})^{T}(P_{\Omega}-\mathbf{x}) \\&~~+(\lambda-\hat{\lambda})^{T}(2 \delta I _{\bar{N}}-K_{n}) \dot{\lambda}
	\\&=2 \delta\partial_{\mathbf{x}}V_{1}^{T}(P_{\Omega}-\hat{\mathbf{x}})-2 \delta\partial_{\mathbf{x}}\hat{V}_{1}^{T}(P_{\Omega}-\hat{\mathbf{x}})\\&~~+2 \delta(\partial_{\mathbf{x}}V_{1}-\partial_{\mathbf{x}}\hat{V}_{1})^{T}(\hat{\mathbf{x}}-\mathbf{x})\\&~~+2 \delta(P_{\Omega}-\mathbf{x})^{T}(P_{\Omega}-\hat{\mathbf{x}})-2 \delta(P_{\Omega}-\mathbf{x})^{T}(P_{\Omega}-\hat{\mathbf{x}})\\&~~+2 \delta(\mathbf{x}-\hat{\mathbf{x}})^{T}(P_{\Omega}-\mathbf{x}) \\&~~+(\lambda-\hat{\lambda})^{T}(2 \delta I _{\bar{N}}-K_{n}) \dot{\lambda}\\&=2 \delta(\partial_{\mathbf{x}}V_{1}^{T}-P_{\Omega}+\mathbf{x})(P_{\Omega}-\hat{\mathbf{x}})-2 \delta\partial_{\mathbf{x}}\hat{V}_{1}^{T}(P_{\Omega}-\hat{\mathbf{x}})\\&~~+2 \delta(\partial_{\mathbf{x}}V_{1}-\partial_{\mathbf{x}}\hat{V}_{1})^{T}(\hat{\mathbf{x}}-\mathbf{x})-2 \delta(P_{\Omega}-\mathbf{x})^{T}(P_{\Omega}-\hat{\mathbf{x}})\\&~~+2 \delta(\mathbf{x}-\hat{\mathbf{x}})^{T}(P_{\Omega}-\mathbf{x})+(\lambda-\hat{\lambda})^{T}(2 \delta I _{\bar{N}}-K_{n}) \dot{\lambda}.
	\end{aligned}\end{equation}
	\vspace{2mm}
	Then use (\ref{pi}), (\ref{cN}) and (\ref{keys}), one can find that $(\partial_{\mathbf{x}}V_{1}^{T}-P_{\Omega}+\mathbf{x})(P_{\Omega}-\hat{\mathbf{x}})\le 0$ and $\partial_{\mathbf{x}}\hat{V}_{1}^{T}(P_{\Omega}-\hat{\mathbf{x}}) \le0$. To procced, caculate
	\begin{equation} \label{la} \begin{aligned}
	\dot{V}_{3}+\dot{V}_{21}	&\le -2 \delta(P_{\Omega}-\mathbf{x})^{T}(P_{\Omega}-\hat{\mathbf{x}})+2 \delta(\mathbf{x}-\hat{\mathbf{x}})^{T}(P_{\Omega}-\mathbf{x}) \\&~~~+2 \delta(\partial_{\mathbf{x}}V_{1}-\partial_{\mathbf{x}}\hat{V}_{1})^{T}(\hat{\mathbf{x}}-\mathbf{x})+(\lambda-\hat{\lambda})^{T}(2 \delta I _{\bar{N}}-K_{n}) \dot{\lambda}
	\\&\le  -2 \delta ||P_{\Omega}-\hat{\mathbf{x}}||^{2}+2 \delta(\partial_{\mathbf{x}}V_{1}-\partial_{\mathbf{x}}\hat{V}_{1})^{T}(\hat{\mathbf{x}}-\mathbf{x}) \\&~~~+(\lambda-\hat{\lambda})^{T}(2 \delta I _{\bar{N}}-K_{n}) \dot{\lambda}
	\\&\le  -2 \delta ||P_{\Omega}-\mathbf{x}||^{2}-2 \delta(\partial F(\mathbf{x})-\partial F(\hat{\mathbf{x}}))^{T} (\mathbf{x}-\hat{\mathbf{x}})\\&~~~-2 \delta(\mathbf{x}-\hat{\mathbf{x}})^{T} K_{n}(\mathbf{x}-\hat{\mathbf{x}})-2 \delta(\lambda-\hat{\lambda})^{T} K_{n}(\mathbf{x}-\hat{\mathbf{x}})\\&~~~-2 \delta(\partial G(\mathbf{x})P_{+}^{\mu}-\partial G(\hat{\mathbf{x}})\hat{\mu})^{T}  (\mathbf{x}-\hat{\mathbf{x}})\\&~~~+(\lambda-\hat{\lambda})^{T}(2 \delta I _{\bar{N}}-K_{n}) \dot{\lambda}\\&\le -2 \delta ||P_{\Omega}-\mathbf{x}||^{2}-2 \delta(\partial G(\mathbf{x})P_{+}^{\mu}-\partial G(\hat{\mathbf{x}})\hat{\mu})^{T}  (\mathbf{x}-\hat{\mathbf{x}})\\&~~~-2 \delta(\mathbf{x}-\hat{\mathbf{x}})^{T} K_{n}(\mathbf{x}-\hat{\mathbf{x}})-(\lambda-\hat{\lambda})^{T}K_{n} \dot{\lambda}.
	\end{aligned}\end{equation}
	
	Then, show the method to tackle the last inequality of (\ref{la}): Under Assumption \ref{a1}, $F(\mathbf{x})$ is convex, then $(\partial F(\mathbf{x})-\partial F(\hat{\mathbf{x}}))^{T} (\mathbf{x}-\hat{\mathbf{x}})\ge 0$. According to $\dot{\lambda}=K_{n}\mathbf{x} $ and Lemma \ref{e}, $K_{n}(\mathbf{x}-\hat{\mathbf{x}})=\dot{\lambda}$. For simplicity, we set $\mathbf{P}_{1}=-2 \delta(\partial G(\mathbf{x})P_{+}^{\mu}-\partial G(\hat{\mathbf{x}})\hat{\mu})^{T} (\mathbf{x}-\hat{\mathbf{x}})$ and $\mathbf{P}_{2}=-2 \delta(\mathbf{x}-\hat{\mathbf{x}})^{T} K_{n}(\mathbf{x}-\hat{\mathbf{x}})-(\lambda-\hat{\lambda})^{T}K_{n} \dot{\lambda}$.
	
	Next, one can calculate:
	\begin{equation}  \begin{aligned}\mathbf{P}_{2}+\dot{V}_{22}&=-2 \delta(\mathbf{x}-\hat{\mathbf{x}})^{T} K_{n}(\mathbf{x}-\hat{\mathbf{x}})-(\lambda-\hat{\lambda})^{T}K_{n} \dot{\lambda}\\&~~~+(\mathbf{x}-\hat{\mathbf{x}})^{T} K_{n}\dot{\lambda}+(\lambda-\hat{\lambda})^{T} K_{n}\dot{\lambda}\\&=-\delta(\mathbf{x}-\hat{\mathbf{x}})^{T} (2\delta K_{n}-K_{n}^{2})(\mathbf{x}-\hat{\mathbf{x}}).\end{aligned}\end{equation}
	\begin{equation}\label{16}   \begin{aligned}	\mathbf{P}_{1}+\dot{V}_{23}+\dot{V}_{4}&=-2 \delta(\partial G(\mathbf{x})P_{+}^{\mu}-\partial G(\hat{\mathbf{x}})\hat{\mu})^{T}  (\mathbf{x}-\hat{\mathbf{x}})\\&~~~+\delta\left(P_{+}^{\mu}-\hat{\mu}\right)^{T}\left(P_{+}^{\mu}-\mu\right)+\delta(\mu-\hat{\mu})^{T}\left(P_{+}^{\mu}-\mu\right)\\&=-2 \delta(\partial G(\mathbf{x})P_{+}^{\mu}-\partial G(\hat{\mathbf{x}})\hat{\mu})^{T}  (\mathbf{x}-\hat{\mathbf{x}})\\&~~~-\delta\left\|P_{+}^{\mu}-\mu\right\|^{2}+2 \delta\left(P_{+}^{\mu}-\hat{\mu}\right)^{T}\left(P_{+}^{\mu}-\mu-G(\mathbf{x})\right) \\&~~~+2 \delta\left(P_{+}^{\mu}-\hat{\mu}\right)^{T}\left(G(\mathbf{x})-G(\hat{\mathbf{x}})\right)
	\\&~~~+2 \delta\left(P_{+}^{\mu}-\hat{\mu}\right)^{T}G(\hat{\mathbf{x}}) . \end{aligned}\end{equation}
	Combining with (\ref{pi}), (\ref{cN}) and (\ref{zong3}), one can conclude that $\delta\left(P_{+}^{\mu}-\hat{\mu}\right)^{T}\left(P_{+}^{\mu}-\mu-G(\mathbf{x})\right)\le0$ and $2 \delta\left(P_{+}^{\mu}-\hat{\mu}\right)^{T}G(\hat{\mathbf{x}}) \le0$. 
	
	To proceed, one can derive that
	\begin{equation}\label{17}  \begin{aligned}
	2\delta&\left(P_{+}^{\mu}-\hat{\mu}\right)^{T}\left(G(\mathbf{x})-G(\hat{\mathbf{x}})\right)\\&\le  2 \delta\left(P_{+}^{\mu}\right)\left[G(x)-G(\bar{x})-\partial G(x)^{T}(x-\hat{x})\right] \\ &~~~+2 \delta\left(P_{+}^{\mu}\right)^{T} \partial G(x)^{T}(x-\hat{x}) \\ &~~~-2 \delta(\hat{\mu})^{T}\left[G(x)-G(\hat{x})-\partial G(\hat{x})^{T}(x-\hat{x})\right]  \\ &~~~-2 \delta(\hat{\mu})^{T} \partial G(\hat{x})^{T}(x-\hat{x}) \\ &\leq  2 \delta\left(P_{+}^{\mu}\right)^{T} \partial G(x)^{T}(x-\hat{x})  -2 \delta(\hat{\mu})^{T} \partial G(\hat{x})^{T}(x-\hat{x}). \end{aligned}\end{equation}
	
	The last inequality holds by the convexity of $g_{i}$ ($i=1,...,N$). Next, combining (\ref{16}) with (\ref{17}) results in
	\begin{equation}\label{18}
	\mathbf{P}_{1}+\dot{V}_{23}+\dot{V}_{4}\le -\delta\left\|P_{+}^{\mu}-\mu\right\|^{2}. \end{equation}
	To sum up (\ref{11})-(\ref{18}), one finally gets :
	\begin{equation} \label{V}	\dot{V}\le -\delta\left\|P_{+}^{\mu}-\mu\right\|^{2}-\delta(\mathbf{x}-\hat{\mathbf{x}})^{T} (2\delta K_{n}-K_{n}^{2})(\mathbf{x}-\hat{\mathbf{x}}).\end{equation}
	
	Consequently, (\ref{V}) is negative semi-definite and $V(t)\le V(0)$. It is easy to get that $V$ is coercive, by a contradiction, 1) holds. Furthermore, using the global invariance principle, we can find that $(\mathbf{x}, \lambda, \mu)$ converges to the largest invariant set in the closure of $S:=\{\mathbf{x} \in V(0): \dot{V}=\mathbf{0}\}$. To analyze set $S$, one can obtain from (\ref{V}) that $\dot{V}=\mathbf{0}$ results in that the point $(\mathbf{x}, \lambda, \mu)$ is actually an equilibrium of (\ref{jbop}). Hence, 2) holds.
\end{IEEEproof}

\section{Simulation results}

In this Section, we consider the following optimization problem to illustrate the efficiency of the main results:
\begin{example}
	
	\begin{equation}
	\begin{aligned}\min \sum_{i=1}^{3}&f_{i}(x_{i}),\\
	~s.t.~&x_{i}^{(v_{1})}=x_{j}^{(v_{1})},\\&g_{i}(x_{i})\le 0, \\& x_{i} \in \Omega_{i},~ i,j=1,2,3.
	\end{aligned}
	\end{equation}
	More specificly, $x_{1}\in \mathbf{R}$, $x_{2}=col[x_{21},x_{22}]\in \mathbf{R}^{2}$ and $x_{3}=col[x_{31},x_{32}]\in \mathbf{R}^{2}$. $\Omega_{1,2,3}=[1,2]$. $\Omega=\{x_{1},x_{2},x_{3}\}$ and $\mathbf{x}=col[x_{1},x_{2},x_{3}]$. $v_{1}=\{1\}$ leads to $x_{1}=x_{21}=x_{31}$.
	$$
	f_{i}(x_{i})=\left\{\begin{array}{lc}
	(x_{i}-1.5)^{2}+|x_{1}-0.5|& i=1,\\
	|x_{i1}-1|+(x_{i2}-1.5)^{2}& i=2,\\|x_{i1}-1|+|x_{i2}-1.5|& i=3.\end{array}\right.
	$$

	$$
	g_{i}(x_{i})=\left\{\begin{array}{cc}
	x_{i}-2& i=1,\\
	e^{x_{i2}}-5& i=2,\\x_{i1}-x_{i2}-0.4& i=3.\end{array}\right.
	$$
	
\end{example}
We select the Laplacian matrix $L_{3}=\left[
\begin{array}{ccc}
-1& 1&0\\
1& -2 &1\\
0&1&-1\\
\end{array}
\right] $. 

\begin{figure}[!hbt]\centering
 	\includegraphics[width=8cm,height=5.6cm]{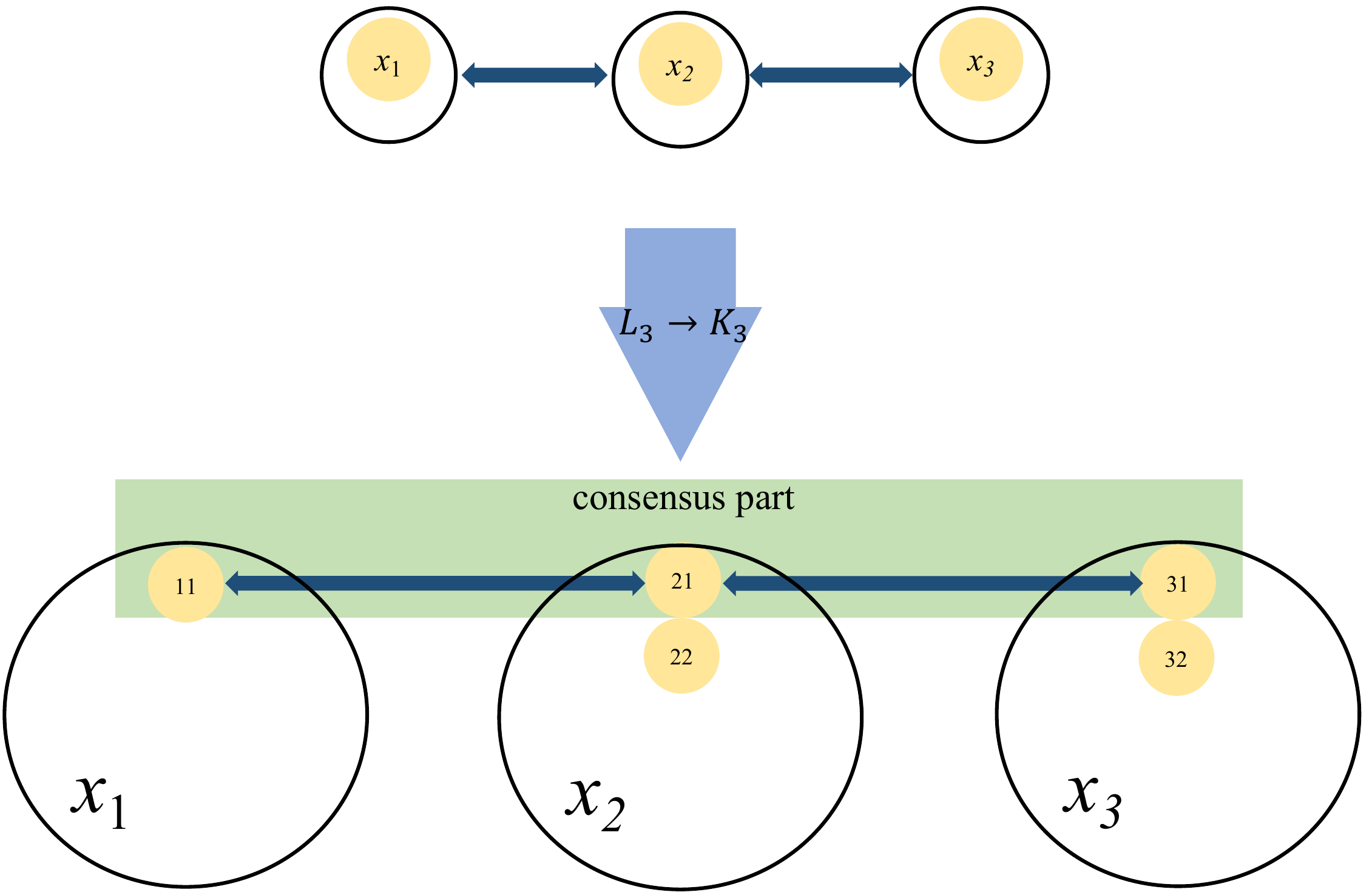}
 	\caption{The graph of corresponding adjacency matrix of $L_{3}$ and $K_{1}$ in Example 2.}
 	\label{to2}
 \end{figure}

 Then $K_{1}$ can be selected as follows:
$$
K_{1}=\left[
\begin{array}{ccccc}
-1 & 1 & 0& 0& 0\\
1& -2& 0& 1& 0 \\
0& 0 & 0 & 0& 0\\
0&1& 0& -1& 0\\
0& 0& 0& 0& 0\\

\end{array}
\right] .
$$
 The graph of $L_{3}$ and $K_{1}$'s corresponding adjacency matrix is shown in Fig. \ref{to2}.
\begin{figure}[!hbt]\centering
	\includegraphics[width=9.5cm,height=7.2cm]{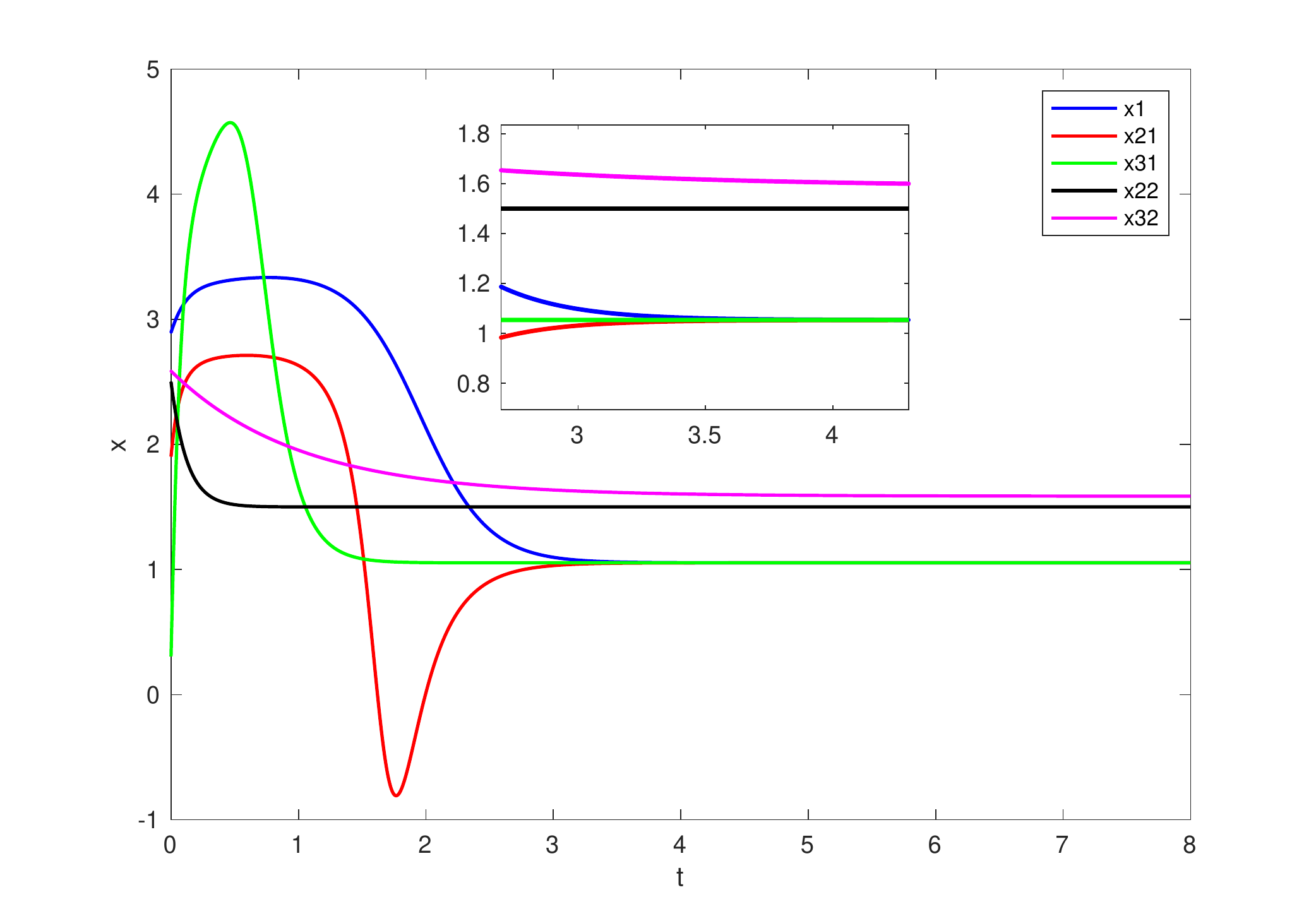}
	\caption{The trajectories of the components of $\mathbf{x}$.}
	\label{x} 	
\end{figure}
\begin{figure}[!hbt]\centering
	\includegraphics[width=9.5cm,height=7.2cm]{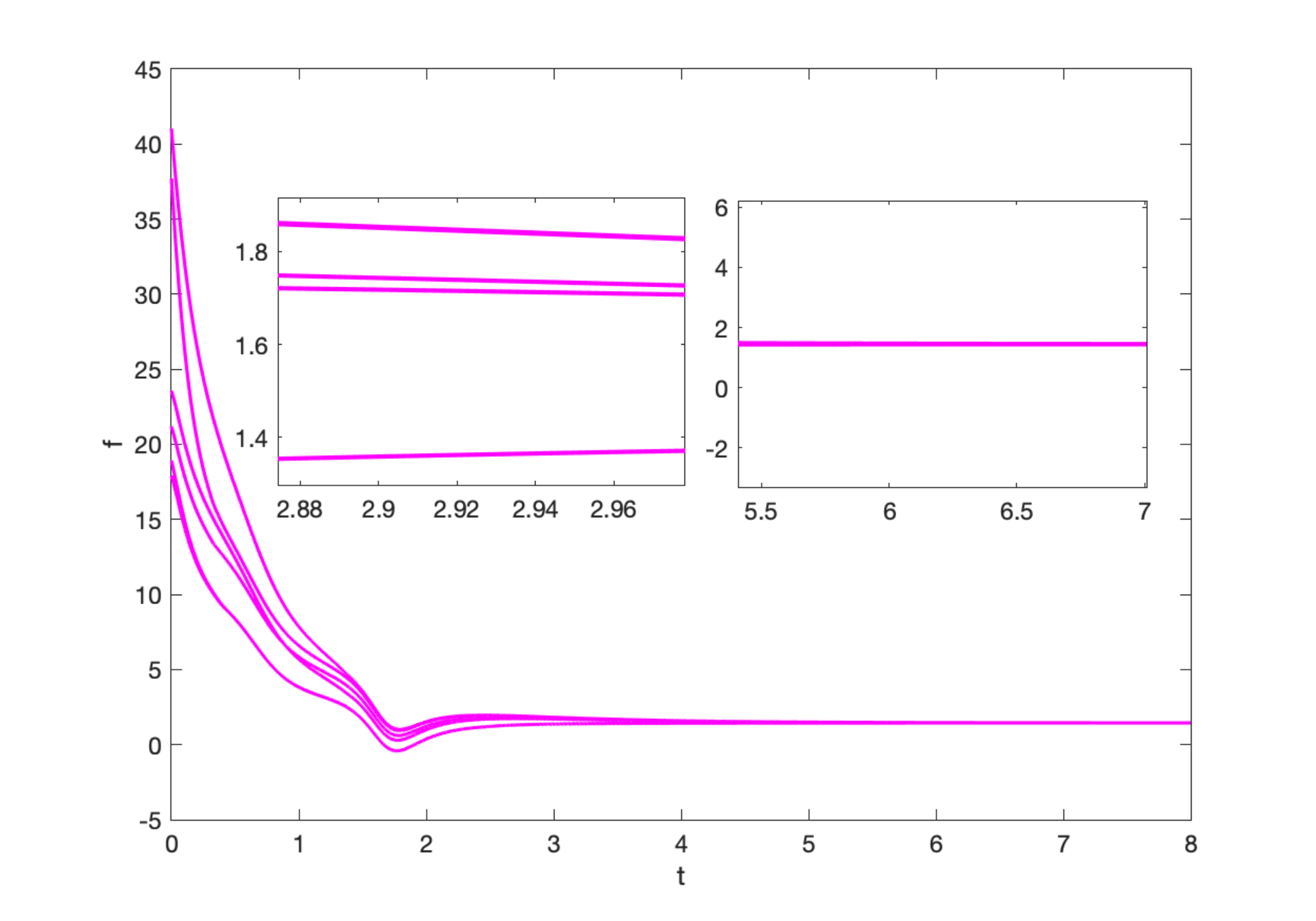}
	\caption{The value of the objective function from different initial states.}
	\label{f} 	
\end{figure}

The solutions of the problem in Example 2 by Algorithm (\ref{maina1}) are shown by Fig.\ref{x} and Fig.\ref{f}. Fig.\ref{x} shares the trajectories of the components of $\mathbf{x}$ and Fig.\ref{f} shares the value of the objective function from different initial states. We can get the solution $x_{1}=1.0536$, $x_{2}=[1.0536, 1.5858]^{T}$ and $x_{3}=[1.0536, 1.5]^{T} $ which can be verified that this point is optimal. The optimal value is 1.4532 from Fig.\ref{f}.

\section{Conclusions}\label{Conclusions}

In this article, a continuous-time decentralized algorithm in RNNs-based fashion has been proposed for the decentralized-partial-consensus optimization under the convex inequality constraints and local set constraints. The partial-consensus matrix was presented to tackle the partial-consensus constraints. By the means of nonsmooth analysis and Lyapunov-based technique, the converge of the designed algorithm has been proved. Finally, a numerical simulation has also been illustrated to show the algorithm performance. In future works, the combination with the matrix graph theory will be considered in more flexible cases for DPCO problems.

	\bibliographystyle{IEEEtran}
	\bibliography{xia.bib}

\end{document}